\documentclass[12pt]{amsart}
\usepackage{amssymb}
\newfont{\ams}{msbm10 at 12pt}
\newfont{\amsi}{msbm8}
\newcommand{\Th}{{\bf Theorem.} }
\newcommand{\Lm}{{\bf Lemma.} }
\newcommand{\Def}{{\bf Definition.} }
\newcommand{\Pf}{{\bf Proof.} }

\newcommand{\Co}{{\bf Corollary.} }
\newcommand{\sq}{$\square $}
\newcommand{\fg}{{\mathfrak g} }
\newcommand{\fu}{{\mathfrak u} }
\newcommand{\fn}{{\mathfrak n} }
\newcommand{\fl}{{\mathfrak l} }
\newcommand{\fb}{{\mathfrak b} }
\newcommand{\fp}{{\mathfrak p} }

\newcommand{\cC}{{\mathcal C} }
\newcommand{\cF}{{\mathcal F} }
\newcommand{\cP}{{\mathcal P} }

\newcommand{\cQ}{{\mathcal Q} }

\newcommand{\cS}{{\mathcal S} }
\newcommand{\cb}{{\bf b} }
\newcommand{\bu}{{\bf u} }
\newcommand{\bp}{{\bf p} }
\newcommand{\bl}{{\bf l} }
\newcommand{\ba}{{\bf a} }
\newcommand{\bU}{{\bf U} }
\newcommand{\bB}{{\bf B} }
\newcommand{\bH}{{\bf H} }
\newcommand{\bP}{{\bf P} }
\newcommand{\bL}{{\bf L} }
\newcommand{\cN}{{\mathcal N} }
\newcommand{\cH}{{\mathcal H} }
\newcommand{\BZ}{{\mathbb Z} }

\newcommand{\BQ}{{\mathbb Q} }
\newcommand{\BO}{{\mathcal O} }
\newcommand{\bO}{{\mathbb O} }
\newcommand{\BC}{{\mathbb C} }
\newcommand{\BN}{{\mathbb N} }

\newcommand{\Res}{\mbox{Res}}
\newcommand{\Ind}{\mbox{Ind}}
\newcommand{\Inf}{\mbox{Inf}}

\begin{document}
\title{Support varieties for quantum groups.}
\author{Viktor Ostrik}
\thanks{The author is partially supported by the U.S. Civilian
 Research and Development Foundation under Award No. RM1-265.}
\address{Independent Moscow University, 11 Bolshoj Vlasjevskij per.,
 Moscow
121002 Russia}
\date{November 1997}
\email{ostrik@nw.math.msu.su}
\maketitle

\section{Introduction.}
Given a root datum $(Y,X,\ldots )$ and an odd integer
$l$ G.Lusztig defined a finite-dimensional Hopf algebra $\bu $
over the cyclotomic field $k:=\BQ (\xi )$ where $\xi$ is a
primitive $l-$th root of unity, see \cite{Lu1}. Suppose that
$l>h$ where $h$ is the Coxeter number of the split reductive
Lie algebra $\fg $ over $k$ associated to the root datum
$(Y,X,\ldots )$. In this case the cohomology $H^{\bullet}(\bu,k)$ of
algebra $\bu $ with trivial coefficients were computed
by V.Ginzburg and S.Kumar, see \cite{GK}. For any $\bu -$module
$M$ one defines the {\em support variety} $|\bu |_M$ (copying the
corresponding constructions from the theory of modular
representation of finite groups and from the theory of
restricted Lie algebras), see \ref{sup} below. It follows from the
result of Ginzburg and Kumar that $|\bu |_M$ is a closed conical
subvariety of the nilpotent cone $\cN \subset \fg$.

Let $\bU_{\xi}$ be the quantum group with divided powers
 associated with
root datum $(Y,X,\ldots )$. Consider the category $\cC $ of
finite-dimensional $\bU_{\xi}-$modules of type {\bf 1}. Any
$M\in \cC$ is a $\bu -$module in natural way. So we can consider
the support variety $|\bu |_M$. It is a $G-$stable subset of
nilpotent cone, where $G$ is a split algebraic group over $k$
 associated
with $(Y,X,\ldots )$. Let $2a$ denote the codimension of $|\bu |_M$
 in $\cN $. Our
first main result (Theorem \ref{dim}) relates the integer $a$
and the dimension of $M$. Namely, we prove that if $l$ is prime then
 $\dim M$
is divisible by $l^a$.

Recall that the category $\cC$ is a direct sum of some its
 subcategories
called linkage classes, see \cite{APW}. Our second main result
 (Theorem \ref{est})
gives an {\em a priori} `estimate' of support variety for a module
 in a fixed linkage class.

As an application we compute (combining Theorems \ref{dim} and
 \ref{est})
the support varieties of Weyl modules, see Theorem \ref{J}. Another
application (and third main result of this note) is the computation
of support varieties for {\em tilting} modules over the quantum
 group
of type $A_n$, see Theorem \ref{til}. This computation verifies the
 remarkable Humphreys'
Conjectures (in quantum setting) on support varieties of tilting
modules and Lusztig's bijection between two-sided cells in affine
Weyl group and nilpotent orbits, see \cite{H} and \ref{map} below.
In fact, the attempts to handle these Conjectures were the starting
point of this work.

Our considerations (except for the computation of
support varieties of tilting modules) are applicable as well to the
situation in characteristic $p>0$. In this context the algebra $\bu$
becomes the {\em restricted enveloping algebra} of the Lie algebra
 $\fg$
(considered over a field of characteristic $p$). In this case our
Theorem \ref{J} gives a positive answer to Jantzen's Question
\cite{J3} 2.7 (1) for $p>h$.

I am happy to express my sincere gratitude to J.Y.Shi who
taught me everything I know about cells in affine Weyl group of
type $A_n$. Moreover, this work arose from the attempts to apply
his results on cells to Humphreys' Conjectures. Thanks are also due
to M.Finkelberg for many helpful conversations and infinite
patience. I am greatly indebted to J.C.Jantzen who found a gap in
the first version of this note. I am grateful to J.Humphreys for
many useful comments.

\section{Preliminaries.}

\subsection{}
We fix a ground field $k$. From now on all objects and tensor
products will be considered over $k$ if not indicated otherwise.

\subsection{}
We fix a simply connected root datum $(Y,X,\ldots )$
corresponding to an indecomposable Cartan matrix.
Let $G,\fg ,\fb ,\fn ,W, R, R_+, \check{R}$ denote the corresponding
split algebraic group, its Lie algebra, its Borel subalgebra,
nilpotent radical of Borel subalgebra, Weyl group, root system, the
set of positive roots,  and the dual root system respectively.
For any $\alpha \in R$ let $\check{\alpha}$ denote the corresponding
 coroot.
Let $\langle \;,\, \rangle :X\times Y\to \BZ$ denote the canonical
 pairing.
Let $\rho \in X$ denote the halfsum of all positive roots.
We will denote by dot (e.g. $w\cdot \lambda$) the action of $W$ on
 $X$ centered in $(-\rho)$.
The set of simple coroots is denoted by $I$.
For any $i\in I$ let $\alpha _i$ denote the corresponding simple
root. Let
$X_+=\{ \lambda \in X|\langle \lambda ,\check{\alpha}_i\rangle \ge 0, i\in I\}$
denote the set of dominant weights.
Let $h$ denote the Coxeter number of $\fg $.

\subsection{}
We fix an odd integer $l>h$, for $\fg$ of type $G_2$ we suppose
that $l$ is not divisible by 3. Let $\xi$ be a primitive $l-$th
root of unity.

Let $\bU_{\xi}$ denote the quantum group with divided powers
corresponding to $(Y,X,\ldots )$
as defined in \cite{LuB} and let $\cC$ denote the category of
finite-dimensional $\bU_{\xi}-$modules of type {\bf 1}.
 Let $E_i^{(s)},F_i^{(s)},K_i^{\pm 1},i\in I,s>0$ denote the
standard  generators of algebra $\bU_{\xi}$. We set
$$
\bU :=\bU_{\xi}/\mbox{Ideal generated by} \{ K_i^l-1\}.
$$
 Let $\bB\subset \bU$ denote the Borel subalgebra.
For any $\lambda \in X, k\ge 0$ let $H^k(\lambda)\in \cC$ denote the
 higher induced module
from one-dimensional representation of $\bB$ corresponding to
 $\lambda$,
see \cite{APW}. For $\lambda \in X_+$ the socle $L(\lambda)$ of
 $H^0(\lambda)$
is simple and nonzero. The modules $L(\lambda)$ are distinct and any
 simple
module in $\cC$ is isomorphic to some $L(\lambda)$.

\subsection{}
Let $W_a$ denote the affine Weyl group: it
is generated by $W$ and by translations by $l\lambda$ for all
 $\lambda \in Y$.
Let $X/W_a$ denote the set of $W_a-$orbits in $X$ with respect to
the dot action. For any $\Omega \in X/W_a$
let $\cC (\Omega)$ denote the full subcategory of $\cC$ consisting
of modules with composition factors $L(\lambda)$ where
$\lambda \in \Omega$.
The Linkage Principle (see \cite{APW}) asserts that
$$
\cC =\bigoplus_{\Omega \in X/W_a}\cC (\Omega).
$$

Moreover all composition factors of any module $H^i(\lambda)$ are
of the form $L(\mu)$ where $\mu \in W_a\cdot \lambda$.

\subsection{}
\label{and}
We refer the reader to \cite{An} for the definition and properties
of tilting modules over $\bU$. Recall that

(i) Any tilting module is a sum of indecomposable ones.

(ii) Indecomposable tilting modules $Q(\lambda)$ are
naturally labelled by the dominant weights $\lambda \in X_+$.

(iii) Tensor product of tilting modules is tilting.

(iv) The module $H^0(\lambda)$ is tilting iff it is simple.

\subsection{}
Let $\cQ$ denote the category of tilting modules over $\bU$.
The full subcategory $\cQ '\subset \cQ$ is called {\em tensor ideal}
if for any $M\in \cQ '$ and $N\in \cQ$

(i) any direct summand of $M$ lies in $\cQ '$.

(ii) $M\otimes N \in \cQ '$.

All tensor ideals in $\cQ$ were described in \cite{Os} using
results of \cite{S1} and \cite{S2} (these results are proved for
simply laced root systems; the proof for non-simply-laced case
is not yet written down).

We will say that weights $\lambda ,\mu \in X_+$ lie in the same
{\em weight cell} if the tensor ideals generated by
$Q(\lambda)$ and $Q(\mu)$ coincide, i.e. $Q(\lambda)$ is a direct
summand of $Q(\mu)\otimes Q$ for some $Q\in \cQ$ and conversely.
 Thus $X_+$ is partitioned
into weight cells (this partition depends on $l$).
It follows from \cite{Os} that the set
of weight cells is bijective to the set of canonical right cells
in an affine Weyl group of $\fg$, see \cite{LuX}.

\subsection{}
Let $\bu $ denote the quantized restricted enveloping
 algebra, see \cite{Lu1}. It is a Hopf algebra of
dimension $l^{\dim \fg}$. Let $\cb$ denote its Borel subalgebra.

\subsection{}
G.Lusztig defined the Frobenius map $\bU \to U(\fg)$,
see \cite{LuB}. Its kernel is $\bU \cdot \bu_+=\bu_+\cdot \bU$
where $\bu_+$ is the augmentation ideal of $\bu$. In particular,
there is a natural $U(\fg)-$module structure on the cohomology
$H^{\bullet}(\bu ,M)$ of a $\bU-$module $M$.

\section{Support varieties.}
In this section we define support varieties
for Hopf algebras (repeating the definitions for
finite groups or for restricted Lie algebras) and
recall their basic properties.

\subsection{}
\label{sup}
Let $\ba$ be a finite dimensional Hopf algebra.
It is known that the cohomology algebra with trivial
coefficients $H^{\bullet}(\ba )=H^{\bullet}(\ba ,k)$ is commutative,
 see \cite{GK}.
Suppose that $H^{\bullet}(\ba )$ is a finitely generated
algebra. Denote by $|\ba|$ the affine algebraic
variety associated to the algebra $H^{ev}(\ba )$ (cohomology of even
 degrees).
The variety $|\ba |$ is equipped with natural $\BC ^*-$
action and distinguished point ${\bf 0}\in |\ba|$.
Let $M$ be a finite dimensional $\ba-$module. Let
 $|\ba|_M\subset |\ba|$ be the zero locus
of the kernel of the natural morphism
$H^{ev}(\ba)\to Ext^{\bullet}_{\ba}(M,M)$.

\Def {\em The variety $|\ba|_M$
is called the support variety of the module $M$.}

It is clear that $|\ba|_M$ is conical with respect to the natural
$\BC ^*-$action and $M\ne 0$ implies that ${\bf 0}\in |\ba|_M$.

{\em Remark.} Let $\{ L_k\}, k\in K$ be a collection
of representatives of all isomorphism classes of
simple $\ba-$modules. Then $|\ba|_M$ is the union of supports
of sheaves $H^{\bullet}(\ba,L_k\otimes M)$ in $|\ba|$. The proof
repeats word by word the arguments in \cite{FP2} 1.4.

\subsubsection{Examples.}
\label{exam}
 (i) Consider the case of the algebra $\ba=\bu$.
In this case it is known that
 $H^{\bullet}(\bu)=H^{ev}(\bu)=\BC [\cN ]$
where $\BC [\cN ]$ denote the algebra of functions on the nilpotent
cone $\cN \subset \fg$, see \cite{GK}. Hence, $|\bu|=\cN $.

(ii) Let $\ba=\cb$. In this case
 $H^{\bullet}(\cb)=H^{ev}(\cb)=S^{\bullet}(\fn ^*)$,
see {\em loc. cit.} Hence, $|\cb |=\fn$.

(iii) Let $M$ be the trivial one-dimensional $\ba-$module. Then
obviously $|\ba|_M=|\ba|$.

(iv) If $M$ is projective then $|\ba|_M={\bf 0}$.

\subsection{}
\label{fun}
The following Lemma is a Hopf-algebra analogue of the
properties of support varieties established in \cite{J1}
in the case of restricted Lie algebras.

\Lm (cf. \cite{J1} 3.3(3), 2.2(1), 3.3(1))
(i) {\em Let $M,N$ be $\ba-$modules.
 Then $|\ba|_{M\oplus N}=|\ba|_M\cup |\ba|_N$.}

(ii) {\em Let $\ba_1, \ba_2$ be Hopf algebras. Let
 $\phi :\ba _1\to \ba _2$ be
a homomorphism of algebras. It defines the morphism
$\phi :|\ba _1|\to |\ba _2|$. For any $\ba _2-$module $M$ we have
$|\ba _1|_M\subset \phi ^{-1}|\ba _2|_M$.}

(iii) {\em Let $0\to M_1\to M_2\to M_3\to 0$ be an exact sequence
of $\ba-$modules. Then for any permutation $(i,j,k)$ of $(1,2,3)$ we
have $|\ba|_{M_i}\subset |\ba|_{M_j}\cup |\ba|_{M_k}$.}

\Pf
(i) is clear from the Remark \ref{sup}.

(ii) is clear from the functoriality of cohomology.

(iii) is clear from the long exact sequence for cohomology. \sq

\subsection{}
\label{equ}
Let $M$ be a $\bU-$module. We will write $|\bu|_M$ for
$|\bu|_{\Res^{\bU}_{\bu}M}$.

\Lm {\em For any $\bU-$module $M$ the variety $|\bu|_M$
is a $G-$invariant subset of the nilpotent cone $\cN$.}

\Pf
See \cite{GK}. $\square$

\subsection{}
\label{ten}

\Lm (cf. \cite{J1} 3.3(5)) {\em Suppose the map
$H^{ev}(\ba)\otimes H^{ev}(\ba)\to H^{ev}(\ba)$ induced by the
coproduct $\delta :\ba\to \ba\otimes \ba$ corresponds to the diagonal
map $|\ba|\to |\ba|\times |\ba|$. Then for any $\ba-$modules $M$ and
$N$ we have $|\ba|_{M\otimes N}\subset |\ba|_M\cap |\ba|_N$.}

\Pf Clear. $\square$

{\em Remarks.} (i) For all the algebras we consider below the
assumption of Lemma is satisfied, and we will omit the checking.

(ii) It can be conjectured that equality holds:
$|\bu|_{M\otimes N}=|\bu|_M\cap |\bu|_N$, cf. \cite{FP1} 2.1(c).

\subsection{}
\Co {\em Suppose $\lambda, \mu \in X_+$ lie in the same
weight cell. Then $|\bu|_{Q(\lambda)}=|\bu|_{Q(\mu)}$.}

\Pf
$Q(\lambda)$ is a direct summand of $Q(\mu)\otimes Q$ for
some $Q\in \cQ$. By Lemmas \ref{fun} (i) and \ref{ten} we obtain
an inclusion $|\bu|_{Q(\lambda)}\subset |\bu|_{Q(\mu)}$.
By symmetry we have the inverse inclusion. $\square$

\subsubsection{}
\label{map1}
\Co {\em There is a well-defined order-preserving map $\cH$
from the set of canonical right cells to the set of closed
$G-$invariant subsets of $\cN$.}

\subsection{}
\label{map}
Now let us formulate the Humphreys' Conjectures. Recall that
the set of canonical right cells is bijective to the set of
all two-sided cells in the affine Weyl group $W_a$, see \cite{LuX}.
Furthermore, G.Lusztig constructed a bijection between the set of
two-sided cells in $W_a$ and the set of nilpotent orbits in $\fg$,
see \cite{Lu2}. Composing these bijections we obtain a bijection
 between the set
of canonical right cells and the set of closed irreducible
$G-$equivariant subsets of $\cN $. The Conjecture 1 below is a
quantum version of Hypothesis \cite{H} 12, and the
Conjecture 2 modulo Conjecture 1 is close to discussion in the
last paragraph of \cite{H} 12.

{\bf Conjecture 1.} {\em The map $\cH$
coincides with Lusztig's bijection.}

Suppose $l$ is prime.

{\bf Conjecture 2.} {\em Suppose $|\bu |_{Q(\lambda)}$ has
codimension
$2a$ in $\cN$. Then $\dim Q(\lambda)$ is divisible by $l^a.$}

{\em Remark.} In fact, J.Humphreys conjectured a stronger
relation between the
dimensions of tilting modules with regular highest weight and
the dimension of their support varieties.

\section{Dimensions of $\bU-$modules.}
\subsection{}
\label{dim}
Let $d_R=\prod_{\alpha \in R_+}\langle \rho ,\check{\alpha}\rangle$
(denominator of the Weyl dimension formula).
This section is devoted to the proof of the following

\Th {\em Let $M\in \cC$. Suppose
the codimension of $|\bu|_M$ in $|\bu|$
is equal to $2a$. Then $d_R\dim M$ is divisible
by $l^a$.}

{\em Remarks.} (i)The codimension of a closed $G-$invariant
subset of $\cN$ is always even since the number of
nilpotent orbits is finite and any orbit is a
simplectic variety.

(ii) If $l$ is prime then we get divisibility of $\dim M$ by $l^a$
(since in this case $d_R$ is not divided by $l>h$).

(iii) The Theorem implies the Conjecture 2 in \ref{map}.

\subsection{}
\Pf
We will write $|\cb|_M$ for $|\cb|_{\Res^{\bU}_{\cb}(M)}$. It is
a subset of $\fn =|\cb|$, see \ref{exam}.

\Lm {\em $\dim |\cb|-\dim |\cb|_M\ge a$.}

\Pf
By Lemma \ref{equ} we have $|\cb|_M\subset \fn \cap |\bu|_M$.
Further, $|\bu|_M$ is a union of a finite number of nilpotent orbits.
Let $\bO$ be a nilpotent orbit. Then $\bO \cap \fn$ is a lagrangian
subvariety in $\bO$, see \cite{GC}. In particular
$\dim (\bO \cap \fn )=\frac{1}{2}\dim {\bO}$.
 Since $\dim \fn =\frac{1}{2}\dim \cN$ we have
$\dim \fn -\dim (\bO \cap \fn)=\frac{1}{2}(\dim \cN -\dim \bO)$. The
result follows. $\square$

\subsection{}
Now we use an interpretation of the dimension of the support variety
as {\em complexity} of the module, see \cite{FP2} 3.1.

\Lm {\em For any $\cb-$module $N$ there exists a projective
resolution $P^{\bullet}\to N$ and a constant $C$ such that
 $\dim (P^i)\le Ci^{\dim
|\cb|_N-1}, i\ge 1$.}

{\bf Proof} repeats word by word \cite{FP2} 3.2. $\square$

\subsection{}
Recall that the algebra $\cb$ has a natural $Y-$grading,
see \cite{LuB}.
We consider the induced $\BZ -$grading such that any simple root has
degree $1$ (principal grading). The module $M$ obviously has a
$\BZ -$grading compatible with the principal grading on $\cb$ since
$M$ is a restriction of
$\bU-$module (this grading comes from the weight grading of $M$).
 We will say that $\cb-$module is graded if it has
$\BZ -$grading compatible with the principal grading of $\cb$.

\subsection{}
For any finite dimensional $\BZ-$graded space
 $N=\oplus _{i\in \BZ} N_i$
let $\dim_tN$ denote its graded dimension
$\dim _tN=\sum_{i\in \BZ}\dim N_it^i$. Set
$p(t)=\prod _{\alpha \in R_+}
\frac{(1-t^{\langle \rho,\check{\alpha}\rangle l})}{(1-t^{\langle
\rho,\check{\alpha}\rangle })}$.
The following result is well-known.

\Lm {\em Let $P$ be an indecomposable projective $\cb-$module. It
admits grading such that $\dim _tP=t^cp(t)$ for some integer $c$.}

We will refer to the number $c$ as to {\em height} of the graded
module $P$.

\subsection{}
\label{res}
\Lm {\em There exists a projective
resolution $P^{\bullet}\to M$ and a constant $C$ such that}

(i) {\em it is graded;}

(ii) {\em indecomposable projectives with height $c$ can occur
with nonzero multiplicity only in $P^i$ with $0\le i\le c$;}

 (iii) {\em $\dim (P^i)\le Ci^{\dim
|\cb|_M-1}, i\ge 1$.}

\Pf Clear. $\square$

\subsection{}
\label{ser}
Let us return to our module $M$. We may (and will) suppose that
 $\dim _tM$ is a polynomial
with nonzero constant term (shifting the grading).

\Lm {\em There is a constant $C$ such that $\dim _tM=p(t)s(t)$
 where $s(t)=\sum _{i=0}^{\infty}s_it^i$ is
a Taylor series with integral coefficients and
$|s_i|\le Ci^{\dim|\cb|_M-1}, i\ge 1$.}

\Pf
Consider the function $s(t)=\frac{\dim_tM}{p(t)}$. It is
a rational function with poles in certain roots of unity of
degree $d_Rl$ since $p(t)=0$ implies $t^{d_Rl}=1$. Let
$s(t)=\sum_{i=0}^{\infty}s_it^i$ be the Taylor expansion of $s(t)$.
Suppose the inequality in Lemma does not hold. For any $1\le k \le d_Rl$
consider the sequence $p_k(i)=s_{k+id_Rl}$. We claim that all these
sequences are polynomials in $i$ for $i$ large enough. Moreover, for
some $k$ the sequence $p_k(i)$ is polynomial of degree
$\ge \dim |\cb|_M$ by our assumption. It follows that
$\sum_{i=1}^j|s_i|$ increases faster than $C'j^{\dim |\cb|_M+1}$.
But this is impossible since $s(t)$ can be obtained by computing
Euler characteristic of the resolution given by Lemma \ref{res}.
The Lemma is proved.
$\square$

\Co {\em The series $s(t)$ is a rational function with poles of
order less than or equal to $\dim |\cb|_M$.}

\subsection{}
Set $d_R(t)=\prod_{\alpha \in R_+}\frac{1-t^{\langle \rho ,
\check{\alpha}\rangle }}{1-t}$.
Note that the polynomial $d_R(t)p(t)$ has zero of order $\dim \fn$
at any nontrivial $l-$th root of unity. It follows from the Lemma
and Corollary \ref{ser} that the polynomial $d_R(t)\dim _tM$ has
zero of order at least $\dim \fn -\dim |\cb|_M\ge a$ at any
nontrivial $l-$th root of unity. Hence the polynomial
$d_R(t)\dim _tM$ with integral coefficients
is divisible by $\left(\frac{t^l-1}{t-1}\right)^a$. Consequently,
$d_R\dim M$ is divisible by $l^a$. The Theorem is proved.
$\square$

\subsection{}
\label{pri}

{\bf Corollary} (of the proof). {\em In notations of Theorem the
 polynomial $\dim_t(M)$ has a
zero of order $\ge a$ at any primitive $l-$th root of unity.
Conversely, let $\dim_t(M)$ have a zero of order $\le a$ at any
 primitive
$l-$th root of unity. Then the codimension of $|\bu|_M$ in $|\bu|$ is
less than or equal to $2a$.}

{\em Example.} Let $M=H^0(\lambda)$ for $\lambda \in X_+$. By the
Weyl character formula (see \cite{APW}) we have
$$
\dim_t(M)=\prod _{\alpha \in R_+}\frac{t^{\langle \lambda+\rho,
\check{\alpha}\rangle }-1}
{t^{\langle \rho,\check{\alpha}\rangle }-1}.
$$
In particular the order of zero at $l-$th primitive root of unity is
equal to the number of positive roots $\alpha$ such that
$\langle \lambda +\rho ,\check{\alpha}\rangle \in l\BZ$.

\section{Estimation of the support variety.}
In this section we estimate the support variety of
a $\bU-$module in a given linkage class. This is a quantum analogue
of Jantzen's Conjecture \cite{J3} 2.7(1).

\subsection{}
For any $\lambda \in X$ consider the following root subsystem of
 $R$:
$$
R_{\lambda}=\{ \alpha \in R|\langle \lambda +\rho ,
\check{\alpha}\rangle \in l\BZ \}.
$$
It is easy to see that if $\lambda $ and $\mu $ lie in the same
linkage class, then $R_{\lambda}$ and $R_{\mu}$ are conjugate
with respect to the Weyl group action. Fix a linkage class $\Omega$
 containing a weight
$\lambda$. Choose $J\subset I$ such that $R_{\lambda}$
is $W-$conjugate to the parabolic root subsystem $R_J$ (this is
possible by \cite{J3} 2.7).
Let $\fp _J\subset \fg$ denote the corresponding parabolic
subalgebra, and let $\fu _J\subset \fp _J$ (resp. $\fl _J$) denote
its nilpotent radical (resp. Levi subalgebra).

\subsection{}
\label{est}
This section is devoted to the proof of the following

\Th {\em The support variety of any module in the
linkage class $\Omega$ is contained in $G\cdot \fu _J$.}

{\em Remark.} $G\cdot \fu _J$ is an irreducible $G-$equivariant
subvariety of $\cN$. Hence, $G\cdot \fu _J$ is the closure of
a nilpotent orbit $\BO $. This orbit is called a Richardson orbit
(corresponding to $J$). It is known that $\dim G\cdot \fu _J=
\dim \BO =2\dim \fu _J$.

Before we start the proof we need some preparations.

\subsection{}
For any $J\subset I$ we consider the parabolic quantum group
$\bP_J$: it is a Hopf subalgebra of $\bU$ generated by generators
$F_i^{(s)}, K_i^{\pm 1}, i\in I, s>0$ and $E_i^{(s)}, i\in J, s>0$
(in notations of \cite{APW} $\bP_J=\bU(J)$). Also consider the
`Levi subgroup' $\bL_J$: it is a Hopf subalgebra of $\bU$ generated
by $K_i^{\pm 1}, i\in I$ and $E_i^{(s)},F_i^{(s)}, i\in J, s>0$.
We have a canonical surjection $\bP_J \to \bL_J$ which sends
$F_i^{(s)},
i\not \in J, s>0$ to zero and other generators into themselves.
Let $\Inf _{\bL_J}^{\bP_J}$ denote the corresponding inflation
functor from $\bL_J-$modules to $\bP_J-$modules. Obviously, this
 functor is exact.

We define the parabolic and Levi subalgebras $\bp _J$
and $\bl _J$ of $\bu$ in the same way. There is the Frobenius
map $\bP_J\to U(\fp_J)$ (resp. $\bL_J\to U(\fl_J)$) with
kernel $\bP_J\cdot (\bp_J)_+=(\bp_J)_+\cdot \bP_J$
(resp. $\bL_J\cdot (\bl_J)_+=(\bl_J)_+\cdot \bL_J$) where
$(\bp_J)_+\subset \bp_J$ (resp. $(\bl_J)_+\subset \bl_J$)
is the augmentation ideal.

\subsection{}
\label{inf}
Let $\bH=\bL_{\varnothing}$ denote the subalgebra of $\bU$ generated
by  $K_i^{\pm 1},i\in I$.
We have natural surjection $\bB=\bP_{\varnothing}\to \bH$.
For any parabolic subalgebra $\bP$ with Levi subalgebra $\bL$
consider the functors $F_1,F_2:$\{ integrable
$\bH-$modules\} $\to$ \{ integrable \bP-modules\},
 $F_1(M)=\Ind _{\bB}^{\bP}(\Inf _{\bH}^{\bB}(M)),
F_2(M)=\Inf _{\bL}^{\bP}(\Ind _{\bL \cap \bB}^{\bL}
(\Inf _{\bH}^{\bL \cap \bB}(M)))$.
Here we use induction functors defined in \cite{APW}.

\Lm {\em The functors $F_1, F_2$ are isomorphic. In particular, for
any integrable $\bH-$module $M$ and $i\ge 0$ we have: }
$$
R^i\Ind _{\bB}^{\bP}(\Inf _{\bH}^{\bB}(M))=\Inf _{\bL}^{\bP}
(R^i\Ind _{\bL \cap \bB}^{\bL}(\Inf _{\bH}^{\bL \cap \bB}(M)))
$$

\Pf
It is enough to compare the right adjoint functors. $\square$

\subsection{}
\label{koh}
\Lm {\em Let $\bp \subset \bu$ be a parabolic subalgebra. Then
$H^{\bullet}(\bp )=H^{ev}(\bp )$ is the algebra of regular functions
on $\cN \cap \fp $ and the natural morphisms $|\bp|\to |\bu|$
and $|\bp|\to |\bl|$ are the natural inclusion
and projection maps respectively.}

{\bf Proof} is the same as in \cite{GK}. Recall it very briefly.
First step is the computation of cohomology of $\cb$.
There exists a $\fb -$equivariant isomorphism of algebras
$H^{2\bullet}(\cb )=S^{\bullet}(\fn ^*),$ see \cite{GK} \S 2.
Further, there are two spectral sequences $'E$ and $''E$,
converging to the same limit, with $E_2-$terms
$$
'E_2^{p,q}=H^p(\bp ,R^q\Ind _{\bB }^{\bP}(k))
$$
and
$$
''E_2^{p,q}=R^p\Ind _{U(\fb )}^{U(\fp )}(H^{q}(\cb ,k)).
$$
Both sequences collapse at $E_2$ and the second one converges to the
algebra $\BC [\fp \cap \cN]$ of regular functions on $\fp \cap \cN$.
The Lemma follows from this easily.
Let us explain the equality $\Ind _{U(\fb )}^{U(\fp )}
(S^{\bullet}(\fn ^*))=\BC [\fp \cap \cN]$ and vanishing of
$R^i\Ind _{U(\fn )}^{U(\fp )}(S^{\bullet}(\fn ^*)), i>0$.
From an exact sequence of $U(\fb)-$modules $0\to (\fn /\fu )^*
\to \fn ^*\to \fu ^*\to 0$
we get that $S^n(\fn ^*)$ has a filtration with quotients of the
form $S^i((\fn /\fu )^*)\otimes S^{n-i}(\fu ^*)$. Note that $\fu ^*$
 is a restriction
of $U(\fp )-$module and $(\fn /\fu )^*$ is an inflation of
 $U(\fb \cap \fl)-$module
$(\fn \cap \fl)^*$. It follows that
$\Ind (S^{\bullet}((\fn /\fu )^*)\otimes S^{\bullet}(\fu ^*))=
\BC [\fl \cap \cN]\otimes S^{\bullet}(\fu ^*)$ and the higher
inductions vanish. Hence, the natural inclusion
 $\BC [\fp \cap \cN]\hookrightarrow
\Ind _{U(\fb )}^{U(\fp )}(S^{\bullet}(\fn ^*))$ induced by the
 restriction map
$\BC [\fp \cap \cN]\to S^{\bullet}(\fn ^*)$ is an
 isomorphism
(since we have an isomorphism of affine algebraic varieties
$\fp \cap \cN \simeq \fu \times (\fl \cap \cN )$) and the
higher inductions vanish. \sq

\subsection{}
\label{ur}
\Lm {\em Let $M$ be a projective $\bl-$module. Consider $M$
as $\bp-$module via canonical surjection $\bp \to \bl$.
Then $|\bp|_M$ is contained in $\fu$ (unipotent
radical of $\fp$).}

\Pf Follows from the naturality of support variety, see Lemma
\ref{fun} (ii). $\square$

\subsection{}
\label{sp}
Let $L$ be a $\bU-$module.

\Lm (cf. \cite{J}, I.6.12) {\em For any $\bP-$module $M$ there
are two spectral sequences converging to the same limit,
with $E_2-$terms: }
$$
'E_2^{p,q}=H^p(\bu ,L\otimes R^q\Ind _{\bP}^{\bU}(M))
$$
and
$$
''E_2^{p,q}=R^p\Ind _{U(\fp )}^{U(\fg )}(H^q(\bp ,L\otimes M)).
$$
{\em Moreover, the algebra $H^{\bullet}(\bu)$ acts on both spectral
 sequences
commuting with all differentials and the actions of
 $H^{\bullet}(\bu)$ on
$'E_{\infty}$ and $''E_{\infty}$ coincide.}

\Pf Consider the following functors $F_1, F_2$ from the category
of integrable $\bP-$modules to integrable $U(\fg)-$modules:
$F_1=(?)^{\bu}\circ \Ind_{\bP}^{\bU}(?),
 F_2=\Ind_{U(\fp)}^{U(\fg)}(?)
\circ (?)^{\bp}$. We claim that these functors are isomorphic.
 Indeed,
the left adjoint functor of $F_1$ is isomorphic to the left
 adjoint functor
of $F_2$, since $\Res_{\bP}^{\bU}(?)\circ \Inf_{U(\fg)}^{\bU}(?)=
\Inf_{U(\fp)}^{\bP}(?)\circ \Res_{U(\fp)}^{U(\fg)}(?)$. It follows
that for any $\bP-$module $N$ there are two spectral sequences $'E$
 and $''E$ converging
to the same limit, with $E_2-$terms:
$$
'E_2^{p,q}=H^p(\bu,R^q\Ind_{\bP}^{\bU}(N))
$$
and
$$
''E_2^{p,q}=R^p\Ind_{U(\fp)}^{U(\fg)}(H^q(\bp,N)).
$$
We set $N=L\otimes M$. Then $R^q\Ind_{\bP}^{\bU}(L\otimes M)=
L\otimes R^q\Ind_{\bP}^{\bU}(M)$
by the generalized tensor identity, see \cite{APW} 2.16 and
\cite{J} I.4.8. Thus, we obtain required sequences.

Now set $N=k$. Then both spectral sequences collapse at $E_2-$terms
and their $E_{\infty}-$terms are equal to $H^{\bullet}(\bu)$
(for $''E$ this follows from
the transitivity of induction and the proof of Lemma \ref{koh}).
The second assertion of the Lemma follows. \sq

\subsection{}
\label{m1}
\Lm {\em Let $\lambda \in X_+$ and $w\in W$. All composition
factors of $H^i(w\cdot \lambda)$ are irreducibles
$L(\mu)$ with $\mu$ linked to $\lambda$ and $\mu \le \lambda$.
Moreover the composition factor $L(\lambda)$ occurs exactly once in
the module $H^{l(w)}(w\cdot \lambda)$, and it does not occur in
$H^i(w\cdot \lambda), i\ne l(w)$.}

\Pf See \cite{APW} and \cite{J} II.6.15-6.16. $\square$

\subsection{}
\label{sind}
Let $\fp$ be a parabolic subalgebra of $\fg$ and let
$S^{\bullet}(\fp ^*)$ be the algebra of functions on $\fp$.
Consider the variety $G\times _P\fp$
(where $P$ is parabolic subgroup of $G$ with Lie algebra $\fp$).
Let $\BC [G\times _P\fp]$
denote the algebra of functions on $G\times _P\fp$.

\Lm {\em There is an isomorphism of algebras}
 $\Ind_{U(\fp )}^{U(\fg)}(S^{\bullet}(\fp ^*))=\BC [G\times _P\fp]$.

\Pf Let $p:G\times _P\fp \to G/P$ denote the natural projection.
 Then the
l.h.s. equals the global sections of the quasicoherent sheaf
 $p_*\BO [G\times _P\fp ]$.
The Lemma follows. \sq

\subsection{}
\label{she}
Consider a graded vector space $M^{\bullet}$ such that

1)$M^{\bullet}$ is a graded finitely generated
 $S^{\bullet}(\fp ^*)-$module,

2)$M^{\bullet}$ is a graded integrable $\fp -$module,

3) structures in 1) and 2) are compatible with respect to the
coadjoint action of $\fp$ on $S^{\bullet}(\fp ^*)$.

Such a datum is equivalent to a datum of a
$P\times \BC ^*-$equivariant
coherent sheaf $\cS (M^{\bullet})$ on $\fp $.
 Let $supp(M^{\bullet})$ denote the support
of this sheaf. For any $i\ge 0$ consider the graded vector space
$N_i^{\bullet}:=R^i\Ind _{U(\fp)}^{U(\fg)}(M^{\bullet})$. It
satisfies all the conditions 1)-3) (for $\fp$ replaced by $\fg$).
Let us explain how this looks geometrically. Let $q:G\times \fp
\to G\times _P\fp$ be the natural projection. There exists a sheaf
$\cS '$ on $G\times _P\fp$ such that $q^*\cS '=\BO [G]
\otimes \cS (M^{\bullet})$.
By Lemma \ref{sind} we have $\cS (N_i^{\bullet})=R^i\phi _*\cS '$
 where $\phi :G\times _P\fp \to \fg,
(g,p)\mapsto gpg^{-1}$.

\Lm {\em We have $supp(N_i^{\bullet})\subset
 G\cdot supp(M^{\bullet})$.}

\Pf Clear. \sq

\subsection{}
{\bf Proof of the Theorem.} It is enough to check the simple modules
$L(\lambda)$. We proceed by induction in $\lambda$ with respect to
standard ordering on $X$ (i.e. $\lambda \ge \mu$ iff
$\lambda -\mu$ is a sum of simple roots with nonnegative
coefficients).  We assume that the Theorem holds for
all modules  $L(\mu)$ where $\mu \in \Omega ,\mu <\lambda$.
Choose $w\in W$ such that $\langle w(\lambda +\rho ),
\check{\alpha}_i\rangle \in l\BN$
for all $i\in J$. Consider the
 $\bP_J-$module $M=\Ind_{\bB}^{\bP_J}(w\cdot \lambda)$. Note that
 our choice of
$w$ gives vanishing of $R^i\Ind_{\bB}^{\bP_J}(w\cdot \lambda)$ for
 $i>0$ (by
Kempf vanishing, see \cite{APW}). In particular,
 $R^i\Ind _{\bP_J}^{\bU}(M)=H^i(w\cdot \lambda)$.
 Further, note that $\Ind_{\bB}^{\bP_J}(w\cdot \lambda)$ is an
 inflation from $\bL_J-$module by Lemma \ref{inf}
 and as $\bL_J-$module it is projective.

 Consider the spectral sequences given by Lemma \ref{sp} for a
 $\bP_J-$module $M$,
 and $L$ running through the set of simple $\bu-$modules (note that
 any simple
$\bu -$module can be lifted to a $\bU-$module). Suppose that
 $|\bu |_{L(\lambda)}$
is not contained in $G\cdot \fu _J$.
Then by Lemmas \ref{m1} and \ref{fun} (iii) for some $L$ the
 $\BC [\cN ]-$
module $H^{\bullet}(\bu ,H^{l(w)}(w\cdot \lambda)\otimes L)$
 has support not
contained in $G\cdot \fu _J$. It follows that the same is true for
$'E_{\infty}$. But this is impossible since
 $'E_{\infty}= ''E_{\infty}$
and by Lemmas \ref{ur}, \ref{she} the support of $''E_{\infty}$
 is contained in
$G\cdot \fu _J$. The Theorem is proved. $\square$

{\em Remark.} The same proof is applicable as well to the modular
 situation with
$p>h$.

\section{Applications.}
\subsection{}
\label{J}
Let $\lambda \in X_+$.

\Th {\em In the notations of \ref{est} we have
$|\bu|_{H^0(\lambda)}=G\cdot \fu _J$.}

\Pf The inclusion $|\bu|_{H^0(\lambda)}\subset G\cdot \fu _J$
is proved in Theorem \ref{est}. On the other hand we know
from Example \ref{pri} and Remark \ref{est} that
$\dim |\bu|_{H^0(\lambda)}\ge \dim G\cdot \fu _J$.
The variety $G\cdot \fu _J$ is irreducible.
The Theorem follows.
 $\square$

{\em Remark.} This gives also the support varieties for Weyl modules
since the support varieties of a module and its dual coincide.

\subsection{}
From now on we consider the case of quantum $SL_{n+1}$, i.e. we
 suppose
that our root datum $(Y,X,\ldots )$ is of type $A_n$. We identify
 $R_+$
with the set of pairs $(j_1,j_2)\in [n+1]\times [n+1], j_1<j_2$
 (where
$[n+1]$ denotes the set of natural numbers $j$ such that
 $1\le j\le n+1$).

\subsection{}
\label{ex}
Let us fix a partition $p =\{ p_1, \ldots \}$ of $n+1$. Let
$p'=\{ p_1',\ldots \}$ denote the dual partition.

\subsubsection{}
\label{seq}
\Lm {\em There exists a decreasing sequence of numbers
 $x_1>x_2>\ldots >x_{n+1}$
such that $x_i-x_j\in \BZ $ and for $\sum_{i=1}^kp_i'<j\le
\sum_{i=1}^{k+1}p_i'$ we have $x_j-x_{j+p_{k+1}'}=l$.}

\Pf Choose any integers $x_1>x_2>\ldots >x_{p_1'}>x_1-l$
 (This is possible since $l>h=n$).
 Now for any $1\le j\le p_1'$ set $x_{j+\sum_{i=1}^kp_i'}=x_j-kl$.
 $\square$

\subsubsection{}
\label{W}
Let us consider a weight $\lambda=(x_1-x_2-1,\ldots ,x_n-x_{n+1}-1)$
(in coordinates corresponding to the basis of
 fundamental weights) where $x_i$ are
the numbers given by \ref{seq}. Obviously, $\lambda \in X_+$.

\Lm {\em The module $Q(\lambda)$ coincides with the induced module
 $H^0(\lambda)$.}

\Pf It is enough to show that module $H^0(\lambda)$ is irreducible
by \ref{and} (iv).
This follows from the quantum analogue of irreducibility criterion
\cite{J}, II.8.21.
  $\square$

\subsubsection{}
\label{dp}
\Lm {\em For $\lambda$ as above $|\bu |_{H^0(\lambda)}=$} \{ the
 closure of orbit
consisting of nilpotent elements with Jordan blocks of sizes
 $(p_1',\ldots )$\} .

\Pf It is easy to see that in this case $R_{\lambda}$ is of type
$A_{p_1-1}\times A_{p_2-1}\times \ldots$. Let $R_J$ be a parabolic
 root
subsystem which is $W-$conjugate to $R_{\lambda}$. By \cite{J3} 2.6
 we have
$G\cdot \fu _J=$ closure of the orbit consisting of nilpotent
 elements with
Jordan blocks of sizes $(p_1',\ldots )$. The result follows from
Theorem \ref{J}. \sq

\subsection{}

In our proof we will use the explicit description of cells
in the affine Weyl group of type $A_n$ given by J.Shi in terms of
admissible sign types. We consider only the case of dominant
admissible sign types, for the general case see \cite{Shi}.

\subsubsection{}
\Def (\cite{Shi}){\em The map $f:R_+\to \{ +,0\} $,
 is called dominant
admissible sign type (dast for short) if for any $\alpha ,\beta ,
\alpha +\beta \in R_+$ the equality $f(\alpha )=+$ yields
$f(\alpha +\beta )=+$.}

Let $\cF$ denote the set of dasts, let $\cP _n$ denote the set of
partitions of $n+1$. Let us define the map $\pi :\cF \to \cP _n$.
 For any $f\in \cF$ we say that a subset $K\subset [n+1]$
is connected if $f(j_1,j_2)=+$ for any $j_1,j_2\in K\; ,j_1<j_2$.
 Now let $p_1$
be the maximal cardinality of a connected subset of $[n+1]$, $p_1+
p_2$ be the maximal cardinality of a subset of $[n+1]$ which is a
disjoint union of two connected subsets and so on (see \cite{Shi}).
Then by C.Greene's Theorem $p_1,\ldots ,p_k$ ($k$ is the largest
integer such that $p_k$ is nonzero) is a partition of $n+1$, i.e.
$p_1\ge p_2\ge \ldots \ge p_k$.
 We set $\pi (f)=(p_1,\ldots ,p_k)\in \cP_n$.

\subsubsection{}
For any $\lambda \in X_+$ we define $f=f(\lambda)\in \cF$ as
follows: for any $\alpha \in R_+$ we set $f(\alpha)=+$ if
 $\langle \lambda +\rho ,\check{\alpha}\rangle \ge l$
and $f(\alpha)=0$ otherwise. We get a map $f:X_+\to \cF$.


\subsubsection{}
Consider the composition $\pi f:X_+\to \cP _n$. The fibers of
 this
map are unions of alcoves (possibly not closed). According to
 \cite{Shi} \S 18.2 and \S 7.2
the partitions of the set of dominant alcoves into fibers of this
map and into weight cells coincide. Moreover the orders on $\cP _n$
and on the set of canonical right cells also coincide. Further,
let $\pi _1$ denote the map from the set of weight cells to $\cP _n$
induced by $\pi f$. Identify $\cP _n$ with the set of nilpotent
orbits in $SL_{n+1}$ as follows: to a partition $(p_1,\ldots )$
we associate the orbit consisting of elements with Jordan blocks of
sizes $(p_1',\ldots)$ (dual partition). It is well known that the
 map $\pi _1$ coincides with Lusztig's
bijection.

\subsubsection{}
{\em Examples.}
(i) For zero weight $\lambda =0$ we have $\pi f(0)=(1,1,\ldots )$
(recall that $l>n$).

(ii) For the weight $\lambda$ from Lemma \ref{W} we have
$\pi f(\lambda)=
(p_1,\ldots ,p_k)$.

\subsection{}
\label{til}
{\bf Theorem.} {\em For the quantum group of type $A_n$
the map $\cH $  (see \ref{map1})
 coincides with Lusztig's bijection.}

{\bf Proof.}
For $\lambda$ of \ref{W} and \ref{dp} we have
 $|\bu |_{Q(\lambda)}=|\bu |_{H^0(\lambda)}=$
the closure of orbit of type $(p_1',\ldots )$. The
 Theorem is proved. $\square $


\begin{thebibliography}{99}

\bibitem[An]{An} H.H.Andersen {\em Tensor products of quantized
tilting modules,} Comm. Math. Phys. {\bf 149} (1992), 149-159.

\bibitem[APW]{APW} H.H.Andersen, P.Polo, K.Wen, {\em Representations
of quantum algebras,} Invent. Math. {\bf 104} (1991), 1-59.

\bibitem[CG]{GC} N.Chriss, V.Ginzburg, {\em Representation theory
and complex geometry,} Birkh\"auser, Boston (1997).

\bibitem[FP1]{FP1} E.M.Friedlander, B.J.Parshall, {\em Support
varieties for restricted Lie algebras,} Invent. Math. {\bf 86}
(1986), 553-562.

\bibitem[FP2]{FP2} E.M.Friedlander, B.J.Parshall, {\em Geometry
of $p-$unipotent Lie algebras,} J. Algebra {\bf 109} (1987), 25-45.

\bibitem[GK]{GK} V.Ginzburg, S.Kumar, {\em Cohomology of quantum
groups at roots of unity,} Duke Math. Journal, vol. {\bf 69}, No. 1
(1993), 179-198.

\bibitem[H]{H} J.E.Humphreys, {\em Comparing modular representations
of semisimple groups and their Lie algebras,} in volume dedicated to
R.E.Block, Modular Interfaces: Modular Lie Algebras, Quantum Groups,
and Lie Superalgebras, ed. V.Chari and I.Penkov,
 International Press, 1997.

\bibitem[J1]{J1} J.C.Jantzen, {\em Kohomologie von $p-$Lie-Algebren
und nilpotent Elemente,} Abh. Math. Sem. Univ. Hamburg {\bf 56}
(1986), 191-219.

\bibitem[J2]{J3} J.C.Jantzen, {\em Support varieties of Weyl
modules,} Bull. London Math. Soc. {\bf 19} (1987), 238-244.

\bibitem[J]{J} J.C.Jantzen, {\em Representations of algebraic
groups,} Academic Press, Orlando (1987).

\bibitem[Lu1]{Lu1} G.Lusztig, {\em Finite dimensional Hopf algebras
arising from quantized universal enveloping  algebras,}
J.Amer.Math.Soc. {\bf 3} (1990), 257-296.

\bibitem[Lu2]{Lu2} G.Lusztig, {\em Cells in affine Weyl group,}
Algebraic groups and related topics, Advanced Studies in Pure
Math., vol. 6, Kinokuniya and North Holland, Tokyo and Amsterdam
(1985); II, J. of Alg. {\bf 109} (1987), 536-548; III, J. Fac.
Sci. Tokyo U., IA {\bf 34} (1987), 223-243; IV, J. Fac. Sci.
Tokyo U., IA {\bf 36} (1989), 297-328.

\bibitem[Lu]{LuB} G.Lusztig, {\em Introduction to quantum groups,}
Birkh\"auser, Boston (1993).

\bibitem[LuX]{LuX} G.Lusztig, N.Xi, {\em Canonical left cells in
affine Weyl group,} Adv. in Math. {\bf 72} (1988), 284-288.

\bibitem[Os]{Os} V.Ostrik, {\em Tensor ideals in category of tilting
modules,} Transformation groups, {\bf 2}, No. 3 (1997), 279-287.

\bibitem[S1]{S1} W.Soergel, {\em Kazhdan-Lusztig-Polynome und eine
Kombinatorik f\"ur Kipp-Moduln,} Electronic Representation Theory
 {\bf 1} (1997), 37-68.

\bibitem[S2]{S2} W.Soergel, {\em Charakterformeln f\"ur Kipp-Moduln
\"uber Kac-Moody-Algebren,} Electronic Representation Theory
{\bf 1} (1997), 115-136.

\bibitem[Shi]{Shi} J.Y.Shi {\em The Kazhdan-Lusztig cells in certain
affine Weyl groups,} Lect. Notes in Math. {\bf 1179},
Springer-Verlag, Berlin, Heidelberg, New York, 1986.

\end{thebibliography}
\end{document}